\documentclass[10pt,leqno]{article}

\usepackage{amsmath,amsfonts,amscd,amssymb,theorem}

\long\def\comment#1\endcomment{}

\makeatletter
\begingroup
\gdef\th@dotted{\normalfont\itshape
  \def\@begintheorem##1##2{%
        \item[\hskip\labelsep \theorem@headerfont ##1\ ##2.]}%
\def\@opargbegintheorem##1##2##3{%
   \item[\hskip\labelsep \theorem@headerfont ##1\ ##2\ (##3).]}}
\endgroup
\makeatother

\theoremstyle{dotted}

\newtheorem{theorem}{Theorem}[section]
\newtheorem{lemma}[theorem]{Lemma}

\newtheorem{corr}[theorem]{Corollary}

\makeatletter
\begingroup
\gdef\th@upshape{\normalfont
  \def\@begintheorem##1##2{%
        \item[\hskip\labelsep \theorem@headerfont ##1\ ##2.]}%
\def\@opargbegintheorem##1##2##3{%
   \item[\hskip\labelsep \theorem@headerfont ##1\ ##2\ (##3).]}}
\endgroup
\makeatother

\theoremstyle{upshape}

\newtheorem{defn}[theorem]{Definition}

\makeatletter
\renewcommand{\subsection}{\@startsection{subsection}{2}{0pt}{-3ex
plus -1ex minus -0.2ex}{-2mm plus -0pt minus
-2pt}{\normalfont\bfseries}} \makeatother

\makeatletter
\@addtoreset{equation}{section}
\makeatother

\newcommand{\cntrct}                
{\hspace{2pt}\raisebox{1pt}{\text{$\lrcorner$}}\hspace{2pt}}

\newcommand{\proof}[1][Proof.]{\smallskip\noindent{\em #1}}
\def\endproof{\hfill\ensuremath{\square}\par\medskip}

\def\eqref#1{\thetag{\ref{#1}}}

\let\latexref=\ref
\def\ref#1{{\normalfont{\latexref{#1}}}}

\newcommand{\calo}{{\cal O}}

\newcommand{\pp}{\mathfrak{p}}

\newcommand{\6}{\partial}

\newcommand{\codim}{\operatorname{\sf codim}}
\newcommand{\cchar}{\operatorname{\sf char}}

\newcommand{\Frac}{\operatorname{\sf Frac}}

\title{Normalization of a Poisson algebra is Poisson}

\author{D. Kaledin\thanks{Partially supported by CRDF grant
RM1-2354-MO-02.}}

\date{{\em Steklov Math Institute}}

\begin{document}

\maketitle


\section*{Introduction}

It is well-known that functions on a smooth symplectic manifold $M$
are equipped with a canonical skew-linear operation called the
Poisson bracket. The bracket is compatible with multiplication in a
certain precise way. Formalizing this structure, one obtains the
notion of a Poisson algebra (see Definition~\ref{poi.def}).

The definition of a Poisson algebra is quite general; among other
things, it involves no assumption of smoothness. Recently there
appeared good reasons to study Poisson algebras in full
generality. In particular, they seem to be quite useful in the study
of the so-called symplectic singularities initiated by A. Beauville
\cite{B}.

However, while non-trivial Poisson structures on smooth manifolds
have been under close scrutiny for fifty years or more, the general
theory is much less developed. It seems that even the simplest facts
are not known, or at least, not easy to find in the existing
literature. 

The goal of the present note is to prove one of these simple facts
-- namely, we prove that the integral closure of a Poisson algebra
is again Poisson (Theorem~\ref{main}). The exposition is essentially
self-contained. We need a couple of preliminary lemmas which are
definitely not new, but not quite standard, either. To spare the
reader the trouble of extensive book search, we have taken the
liberty of re-proving them from scratch.

\section{Statements and definitions.}

Fix once and for all a base field $k$ of characteristic $\cchar k =
0$.

\begin{defn}\label{poi.def}
A {\em Poisson algebra} over the field $k$ is a commutative algebra
$A$ over $k$ equipped with an additional skew-linear operation
$\{-,-\}:A \otimes A \to A$ such that
\begin{equation}\label{poi}
\{a,bc\} = \{a,b\}c + \{a,c\}b\quad, \quad
0 = \{a,\{b,c\}\} + \{b,\{c,a\}\} + \{c,\{a,b\}\},
\end{equation}
for all $a,b,c \in A$. An ideal $I \subset A$ is called a {\em
Poisson ideal} if $\{i,a\} \in I$ for any $i \in I$, $a \in A$.
\end{defn}

Additionally, we will always assume that a Poisson algebra $A$ has a
unit element $1 \in A$ such that $\{1,a\}=0$ for every $a \in A$.

\begin{defn}
A {\em Poisson scheme} over $k$ is a scheme $X$ over $k$ equipped
with a skew-linear bracket in the structure sheaf $\calo_X$
satisfying \eqref{poi}.
\end{defn}

\begin{lemma}\label{loc}
Let $A$ be a Poisson algebra.
\begin{enumerate}
\item For any multiplicative system $S \subset A$, the localization
$A[S^{-1}]$ carries a canonical Poisson algebra structure.
\item Any associated prime ideal $\pp \subset A$ is a Poisson ideal.
\item The radical $J \subset A$ of the algebra $A$ is a Poisson
ideal.
\end{enumerate}
\end{lemma}

\proof{} For \thetag{i}, set
\begin{align*}
\{a_1s_1^{-1},a_2s_2^{-1}\} &= \{a_1,a_2\}(s_1s_2)^{-1} -
\{a_1,s_2\}a_2\left(s_1s_2^2\right)^{-1}\\
&\quad - \{s_1,a_2\}a_1\left(s_1^2s_2\right)^{-1} +
a_1a_2\{s_1,s_2\}\left(s_1^2s_2^2\right)^{-1}.
\end{align*}
For \thetag{ii}, note that $\pp \subset A$ is the kernel of the
canonical Poisson map from the Poisson algebra $A$ to the fraction
field $A_\pp$. For \thetag{iii}, note that $J$ is the intersection
of all the associated primes.
\endproof

Lemma~\ref{loc}~\thetag{i}, in particular, means that the spectrum
of a Poisson algebra is a Poisson scheme. We also note the following
geometric corollary.

\begin{corr}
Let $X$ be a Poisson scheme over $k$. Then the reduction $X_{red}$
of the scheme $X$ is a Poisson scheme, and so is every irreducible
component $X_0$ of the reduction $X_{red}$.\endproof
\end{corr}

Our main result is the following.

\begin{theorem}\label{main}
Let $A_0$ be a Noetherian domain over $k$, and let $A$ be its
integral closure in its fraction field.
\begin{enumerate}
\item Every derivation $\xi$ of the algebra $A_0$ extends to a
derivation of the algebra $A$.
\item Every Poisson bracket $\{-,-\}$ on the algebra $A_0$ extends
to a Poisson bracket on the algebra $A$.
\end{enumerate}
\end{theorem}

Note that both derivations and Poisson brackets extend naturally and
uni\-que\-ly to the fraction field $\Frac A_0 = \Frac A$. The point
is that both preserve the integral closure $A \subset \Frac A$. The
first claim is well-known; nevertheless, we will prove it, because
it is needed in the proof of \thetag{ii}.

The geometric corollary (in fact, an equivalent geometric
formulation) of Theorem~\ref{main} is the following.

\begin{corr}\label{geom}
Let $X_0$ be a Noetherian integral scheme over $k$, and let $X$ be
its normalization. Then every vector field $\xi$ on $X_0$ and every
Poisson scheme structure on $X_0$ extend to $X$.
\end{corr}

\section{Discrete valuation rings.}

To prove Theorem~\ref{main}, we first study the situation in
codimension $1$. In this section, assume given a local Noetherian
algebra $A_0$ over $k$ of dimension $1$. Let $K_0$ be its residue
field. Let $A$ be the integral closure of the algebra $A_0$. It is
well-known that $A$ is a discrete valuation ring, whose residue
field $K$ is a finite extension of the residue field $K_0$. Denote
the valuation by $v$. Fix a uniformizing element $\pi \in A$,
$v(\pi) = 1$.

\begin{lemma}\label{gen}
There exists a single element $x \in A$ generating $A$ over $A_0$.
\end{lemma}

\proof{} By the Primitive Element Theorem, the field $K$ is
generated over $K_0$ by a single element, say $\overline{x}$. Let
$P(x)$ be the minimal polynomial for $\overline{x}$ over $K_0$.
Lift $\overline{x}$ to an element $x \in A$ and consider
$$
y = P(x) \in A.
$$
By definition, we have $y=0 \mod \pi$, so that $v(y) > 0$. If $v(y)
= 1$, we are done: $x$ and $y$ generate $A$ over $A_0$, and $y =
P(x)$. If not, replace $y$ with
$$
y' = P(x + \pi).
$$
By the binomial formula, we have
$$
y' = P'(\overline{x})\pi \mod \pi^2.
$$
Since the polynomial $P$ is minimal, its derivative $P'$ satisfies
$P'(\overline{x}) \neq 0$. Therefore $v(y') = 1$, and we are done:
$A$ is generated over $A_0$ by $x + \pi$.
\endproof

\begin{lemma}\label{vf}
Every derivation $\xi_0:A_0 \to A_0$ of the algebra $A_0$ extends to a
derivation of the algebra $A$.
\end{lemma}

\proof{} Consider the formal power series algebra $B = A[[t]]$ in
one indeterminate $t$. It is a regular local algebra, in particular,
it is integrally closed (see, for example, \cite[Prop. 14]{AC}).
Therefore it coincides with the integral closure of the power series
algebra $B_0=A_0[[t]]$. By functoriality, every automorphism of the
algebra $B_0$ extends to an automorphism of its integral closure
$B$. Consider the automorphism $\sigma_0:B_0 \to B_0$ given by
$$
\sigma_0(t) = t \qquad \sigma_0(a) = \exp(t\xi)(a) \text{ for }a \in
A_0 \subset B_0.
$$
Extend it to an automorphism $\sigma:B \to B$ of the algebra
$B$. Setting
$$
\xi(a) = \frac{\6}{\6 t}\sigma(a) \mod t
$$
gives a derivation $\xi:A \to A=B/tB$ extending the given derivation
$\xi_0$.
\endproof

\begin{lemma}\label{nrm}
Assume that the algebra $A_0$ is equipped with a Poisson bracket
$\{-,-\}$. Then this bracket extends uniquely to the algebra $A$.
\end{lemma}

\proof{} Extend the bracket to the fraction field $\Frac A$. We have
to prove that $\{f,g\} \in A$ for every $f,g \in A$. Let $x \in A$
be the generator provided by Lemma~\ref{gen}. It suffices to prove
that $\{x,x\} \in A$ and $\{x,f\} \in A$ for every $f \in A_0$. But
$\{x,x\} = 0$ tautologically, and $\{x,f\} \in A$ by Lemma~\ref{vf}
(define $\xi_0:A_0 \to A_0$ by $\xi_0(a) = \{a,f\}$).
\endproof

\section{Proof of the Theorem.}

We can now prove Theorem~\ref{main}. It is more convenient to
approach it in the geometric form of the Corollary~\ref{geom}. Thus,
let $X_0$ be a Noetherian integral scheme, and let $X$ be its
normalization. By Lemma~\ref{vf} and Lemma~\ref{nrm},
Corollary~\ref{geom} holds for the open complement $U \subset X$ to
a subscheme $Z \subset X$ of codimension $\codim Z \geq
2$. Therefore we have a derivation and/or a Poisson bracket on the
structure sheaf $\calo_U$. This induces a derivation and/or a
Poisson bracket on the sheaf $j_*\calo_U$, where $j:U
\hookrightarrow X$ is the embedding. Since $\codim Z \geq 2$, and
$X$ is normal, we have $\calo_X \cong j_*\calo_U$.\endproof

\subsection*{Acknowledgements.} I would like to thank
E. Amerik, R. Bezrukavnikov, A. Kuznetsov, M. Rovinsky and
M. Verbitsky for stimulating discussions.

\bigskip

\noindent
{\sc Steklov Math Institute\\
Moscow, USSR}\\
{\em E-mail\/}: {\tt kaledin@mccme.ru}

\end{document}